\documentclass{ws-ijbc}
\usepackage{ws-rotating}     % used only when sideways tables/figures are used
\begin{document}

\catchline{}{}{}{}{} % Publisher's Area please ignore

\markboth{Jian Ren, Chujin Li,  Xingye Kan, Jinqiao Duan and Ting
Gao}{Mean Exit Time and Escape Probability for a Tumor Growth System
under Non-Gaussian Noise}

\title{Mean Exit Time and Escape Probability for a Tumor Growth \\System under
Non-Gaussian Noise\footnote{This research was partly supported by
the NSF Grants 0731201 and 1025422,  the NSFC grants 10971225 , 11028102 and
10901065, and the Fundamental Research Funds for the Central
Universities, HUST 2010ZD037. }}

\author{Jian Ren }

\address{School of Mathematics and Statistics,
Huazhong University of Sciences and Technology\\ Wuhan 430074, China
\\renjian0371@gmail.com }

\author{Chujin Li}
\address{School of Mathematics and Statistics,
Huazhong University of Sciences and Technology\\ Wuhan 430074,
China\\lichujin@mail.hust.edu.cn}

\author{Ting Gao}
\address{Department of Applied Mathematics,
Illinois Institute of Technology\\ Chicago, IL 60616, USA}

\author{Xingye Kan}
\address{Department of Applied Mathematics,
Illinois Institute of Technology\\ Chicago, IL 60616, USA}

\author{Jinqiao Duan}
\address{Department of Applied Mathematics,
Illinois Institute of Technology\\ Chicago, IL 60616, USA
\\duan@iit.edu}

\maketitle

\begin{history}
\received{(to be inserted by publisher)}
\end{history}

\begin{abstract}
Effects of non-Gaussian $\alpha-$stable L\'evy noise on the Gompertz
tumor growth model are quantified by considering the mean exit time
and escape probability of the cancer cell density from inside a safe
or benign domain. The mean exit time and escape probability problems
are formulated in a differential-integral equation with a fractional
Laplacian operator. Numerical simulations are conducted to evaluate
how the mean exit time and escape probability vary or bifurcates
when $\alpha$ changes.   Some bifurcation phenomena are observed and
their impacts are discussed.
\end{abstract}


%\begin{multicols}{2}
\section{Introduction}
\noindent Mathematical models have been proposed to describe the
evolution of numerous biological phenomena such as cancer growth.
The deterministic Gompertz growth model for a tumor density growth
is given by the following differential equation
\begin{equation}\label{deterministic}
\frac{dx}{dt}=(A-B \ln x)x,
\end{equation}
where $x(t)$ describes the density of cancer cells at time $t$,
parameters $A$ and $B$ denote growth and decay rates respectively.

This is a deterministic model, nevertheless. As there is a
discrepancy between clinical data and basic theory due to subtle
environmental fluctuations, we need a stochastic part to describe
the fluctuation. A stochastic model proposed as a functional
Fokker-Planck equation in consideration of both fission and
mortality has been considered in \cite{Lo07} and a diffusion process
with modified infinitesimal mean has been studied in
\cite{Gutierrez-jaimez07}. In \cite{Lo09} and \cite{Albano09} a
multiplicative noise of Brownian motion has been considered and in
\cite{Albano06} a deterministic parameter has been modified by a
group of random variables to depict the stochastic fluctuation.
\cite{Lo09} compared the effect of different therapies by giving
their probability density functions(pdfs) and \cite{Albano09},
\cite{Albano06} both deal with the first exit time problem.

In the present paper we consider the Gompertz growth model driven by
the non-Gaussian $\alpha-$stable noise with jumps and quantify the
dynamics in terms of mean exit time and escape probability.
Specifically, we consider the evolution of the tumor growth in the
density range that cannot be diagnosed, and compute the likelihood
and mean time for a tumor to become diagnosable or malignant. We
examine how this depends on $\alpha$ values.   To this end, we
investigate the differential equation driven by symmetric
$\alpha-$stable L\'evy process $L_t^\alpha$ with characteristics
$(0, 1, \nu)$ with jump measure $\nu(dy)=dy/|y|^{1+\alpha}$
\begin{equation}
dX=[(A -B \ln X)X]dt+   dL_t^\alpha , \quad\quad X(0)=x.
\end{equation}
where $x$ is the initial tumor density.

 \par The paper is arranged as follows. We present a brief introduction about
 L\'evy process in section 2 and then a stochastic tumor growth
model in section 3. Numerical experiments are conducted in section
4.

\section{L\'evy process}
A L\'evy process $(L_t$, $t\geq 0)$ is a stochastic process defined
on a probability space $(\Omega, \mathcal{F}, P)$ with the following
properties:
\begin{romanlist}[(iv)]
 \item  $L_0=0$ $(a. s.)$;

 \smallskip

 \item  $L_t$ has independent and stationary increments, i.e. for each $n$ and
$0\leq t_1\leq t_2\leq \cdots \leq t_{n+1}<\infty$ the random
variables $(L_{t_{j+1}}-L_{t_j}, 1\leq j\leq n)$ are independent and
each $L_{t_{j+1}}-L_{t_j}$ and $L_{t_{j+1}-t_j}$ are equal in law
both satisfy a non-Gaussian distribution.

\smallskip

 \item   $L_t$ has stochastically continuous, i.e. for all $c>0, s>0$
$$\lim_{t\rightarrow s}P(|L_t-L_s|>c)=0.$$

\smallskip

 \item  The paths of $L_t$ are
$\mathbb{P}-$almost surely right continuous with left limits.
\end{romanlist}

 For the characteristic exponent
$\Psi_t(\lambda)=\log\mathbb{E}(e^{i\lambda L_t})$, with the fact
that $L_t$ has independent and stationary increment we get
$\mathbb{E}(e^{i\lambda
L_t})=e^{t\Psi_1(\lambda)}=:e^{t\Psi(\lambda)}$, therefore in the
next paragraph we refer to $\Psi(\lambda)$ as a characteristic
exponent of the L\'evy process.

\par {\bf L\'evy-Khintchine formula \cite{Applebaum09, Kyprianou06} for
L\'evy process:} Suppose that $a\in\mathbb{R}$, $\sigma\geq0$ and
$\nu$ is a measure concentrated on $\mathbb{R} \setminus\{0\}$ such
that $\int_{\mathbb{R}}(1\wedge x^2)\,\mu(dx)<\infty$. From this
triple $(a, \sigma, \nu)$ defined for each $\lambda\in\mathbb{R}$,
$$\Psi(\lambda)=ia\lambda-\frac{1}{2}\sigma^2\lambda^2-\int_{\mathbb{R}}(1-e^{i\lambda x}
+i\lambda xI_{|x|<1})\,\nu(dx).$$ Then there exits a probability
space $(\Omega, \mathcal{F}, \mathbb{P})$ on which a L\'evy process
is defined having characteristic exponent $\Psi$.
\par It can be seen that after some reorganization the characteristic exponent
 of L\'evy process can be rewritten as
 \begin{eqnarray*}
 \Psi(\lambda)&=&\Big\{ia\lambda-\frac{1}{2}\sigma^2\lambda^2\Big\}\\{}&&-\Big\{\nu(\mathbb{R}
\setminus(-1, 1))\int_{|x|\geq1}(1-e^{i\lambda
x})\,\frac{\nu(dx)}{\nu(\mathbb{R} \setminus(-1,
1))}\Big\}\\{}&&-\Big\{\int_{0<|x|<1}(1-e^{i\lambda x} +i\lambda
x)\,\nu(dx)\Big\}\\{}&=:&\Psi^{(1)}-\Psi^{(2)}-\Psi^{(3)}
 \end{eqnarray*}
 where $\Psi^{(1)}$, $\Psi^{(2)}$ are respectively the characteristic exponent of
 $X_t^{(1)}=-\sigma B_t+at,\quad t\geq0$
and $X_t^{(2)}=\sum_{i=1}^{N_t}\xi_i, \quad t\geq0$, here $\{N_t:
t\geq0\}$ is a Poisson process with rate
$\nu(\mathbb{R}\setminus(-1, 1))$ and $\{\xi_i: t\geq0\}$ are i.i.d
with distribution $\nu(dx)/\nu(\mathbb{R}\setminus(-1, 1))$
concentrated on $\{x: |x|\geq1\}$ (if $\nu(\mathbb{R}\setminus(-1,
1))=0$ process $X_t^{(2)}$ is identically
zero).\\
We have a result from \cite{Kyprianou06} to say the composition of a
L\'evy process.

\par {\bf L\'evy-It\^{o} decomposition:} Given any $a\in
\mathbb{R}$, $\sigma\geq0$ and measure $\nu$ concentrated on
$\mathbb{R}\setminus\{0\}$ satisfying
$\int_{\mathbb{R}\setminus\{0\}}(|x|^2\wedge1)\,\nu(dx)<\infty,$
there exists a probability space on which three independent L\'evy
processes $X^{(1)}$, $X^{(2)}$, $X^{(3)}$ exit, where $X^{(1)}$ is a
linear Brownian motion with drift with characteristic exponent
$\Psi^{(1)}$, $X^{(2)}$ is a compound Poisson process with
characteristic exponent $\Psi^{(2)}$ and $X^{(3)}$ is a square
integrable martingale with an almost surely countable number of
jumps on each finite time interval which are of magnitude less than
unity and with characteristic exponent $\Psi^{(3)}$. By taking
$X=X^{(1)}-X^{(2)}-X^{(3)}$ we have that there exist a probability
space on which a L\'{e}vy process is defined with characteristic
exponent
$$\Psi(\lambda)=ia\lambda-\frac{1}{2}\sigma^2\lambda^2-\int_{\mathbb{R}}(1-e^{i\lambda
x}+i\lambda xI_{(|x|<1)})\,\mu(dx)$$ for $\lambda\in\mathbb{R}$.\\

 A random variable $X$ is said to be stable if for all $n>0$,
 the distribution of $X$ and
$S_n=\frac{X_1+\cdots+X_n-a_n}{\sigma_n}$ are equivalent, where
$X_1, \cdots, X_n$ are a sequence of i.i.d random variable, $a_n\in
\mathbb{R}$ and $\sigma_n>0$ for $n\in\mathbb{N}$.  It is well known
from \cite{Applebaum09, Feller71} that the only possible choice is
$\sigma_n=n^{1/\alpha}$ for $\alpha\in(0, 2]$ so we refer to
$\alpha$ as the index. \\
For $\alpha\in(0, 1)\cup(1, 2)$ stable variables have characteristic
exponents of the form
$$\Psi(\lambda)=ia\lambda-c^\alpha|\lambda|^\alpha(1-i\beta
\tan\frac{\pi\alpha}{2}sgn\lambda),$$ where $\beta\in[-1, 1]$,
$a\in\mathbb{R}$, $c>0$ and the sign function
$sgn\lambda=I_{(\lambda>0)}-I_{(\lambda<0)}$;\\
for $\alpha=1$, the characteristic exponents have the form
$$\Psi(\lambda)=ia\lambda-c|\lambda|(1+i\beta
\frac{2}{\pi}sgn\lambda\log|\lambda|),$$ where
$\beta\in[-1, 1]$, $a\in\mathbb{R}$ and $c>0$;\\
for the special case $\alpha=2$,
$$\Psi(\lambda)=ia\lambda-\frac{1}{2}\sigma^2\lambda^2;$$
and the measure
$$\nu(dx)=
\begin{cases}
       \frac{c_1}{x^{1+\alpha}}dx &for\quad x\in(0, \infty);\\
       \frac{c_2}{(-x)^{1+\alpha}}dx  &for\quad x\in(-\infty, 0),
\end{cases}$$
where  $c=c_1+c_2$, $c_1$, $c_2\geq0$ and
$\beta=\frac{c_1-c_2}{c_1+c_2}$
  if $\alpha\in(0, 1)\cup(1, 2)$.  The case $\beta=0$, i.e. $c_1=c_2$, is called
  symmetric $\alpha-$stable L\'evy  process.\\

\par The (infinitesimal) generator $A$ of a L\'evy process $X(t)$ is defined as $A\varphi=\lim\limits_{t\downarrow 0}
\frac{T_t\varphi-\varphi}{t}$ where
$(T_t\varphi)(x)=\mathbb{E}_x(\varphi(X(t))$ for any $\varphi\in
D_A$(the domain of generator $A$). It has the expression
\cite{Applebaum09}
$$(A\varphi)(x)=a\varphi'(x)+\frac{1}{2}\sigma^2\varphi''(x)+\int_{\mathbb{R}
\setminus\{0\}}[\varphi(x+y)-\varphi(x)-y\varphi'(x)I_{|y|<1}]\,\nu(dy).
$$

\section{Tumor growth under $\alpha-$stable L\'evy noise}

In this section we consider the stochastic model driven by symmetric
$\alpha-$stable L\'evy process, i.e. a L\'evy process $L_t$, in
which each $L_t$ is a stable random variable and here for simplicity
we take $c_1=c_2=1$. We refer to the monoclonal benign parathyroid
tumor and take $A=6.46$ year$^{-1}$, $B=0.314$ year$^{-1}$,
corresponding to a parathyroid tumor mean age of $19.6$ years. Here
we investigate in the range of the tumor density that cannot be
diagnosed, we take the two boundary points $S_1=1$ and
$S_2=1.074\times10^8$ (the tumor cell density $x=1.074\times10^8$ is
the smallest diagnoseable mass), how the tumor density progresses.
When a tumor's density is bigger than $1.074\times10^8$, we say it
becomes malignant, while if its density is smaller than $1$, it is
not diagnosable.

\par Based on the deterministic Gompertz model, the corresponding stochastic model is given
by
\begin{equation}\label{levy model}
dX=[(A -B \ln X)X]dt+ dL_t^\alpha , \quad\quad X(0)=x,
\end{equation}
where $L_t^\alpha$ is with the characteristics $(0, 1, \nu)$. Here
the jump measure   $\nu(dy)=dy/|y|^{1+\alpha}$. It is known
\cite{Applebaum09} that  the infinitesimal generator $A$ for
$L_t^\alpha$ is
\begin{equation}
(Af)(x)=\frac{1}{2}f''(x)+\int_{\mathbb{R}
\setminus\{0\}}\frac{f(x+y)-f(x)-yf'(x)I_{|y|<1}}{|y|^{1+\alpha}}\,dy.
\end{equation}
Due to the L\'evy-It\^{o} decomposition the generator for the
process $X(t)$ in \eqref{levy model} is then
\begin{eqnarray}
Af=(A -B \ln x)xf'(x)+\frac{1}{2}f''(x)+\int_{\mathbb{R}
\setminus\{0\}}\frac{f(x+y)-f(x)-yf'(x)I_{|y|<1}}{|y|^{1+\alpha}}\,dy.
\end{eqnarray}

In this symmetric $\alpha-$stable L\'evy motion case,
  the integral operator is related to the
fractional Laplacian operator \cite{Albeverrio00, Guan05, Chen10,
Caffarelli07}. Namely, for $\alpha \in (0, 2)$
\begin{equation}\label{flaplace}
   \int_{\mathbb{R} \setminus\{0\}} [u(x+  y)-u(x) ] \; \nu_{\alpha}(dy) = c_{\alpha} \; (-\Delta)^{\frac{\alpha}2} u(x),
\end{equation}
where
$$
c_{\alpha} \triangleq  \int_{\mathbb{R} \setminus\{0\}} (\cos y -1 )
\; \nu_{\alpha}(dy) <0.
$$
Note that when $\alpha=2$, $L_2 = \Delta$.

\medskip

We consider two issues: Mean exit time and escape probability. We
quantify the stochastic dynamics using these two tools. Especially,
we examine possible bifurcation phenomena when the stability
parameter $\alpha$ various in $(0, 2)$.

  The mean exit time $$u(x)=\mathbb{E} \,\inf\{t\geq0,
X_t(\omega, x)\notin D\},$$ from a bounded interval $D$, is the
first time that the tumor density $X(t)$, initially at value $x$ in
$D$,  gets out of $D$. The mean exit time $u(x)$, from the
diagnoseable range $D= (1, 1. 074\times10^8)$,  satisfies an
integral-differential equation  \cite{Brannan99}
\begin{eqnarray}
Au(x)&=&-1,\quad x\in D=(1, 1. 074\times10^8), \nonumber\\
u(x)&=&0, \quad x\notin D.
\end{eqnarray}
i.e.
\begin{eqnarray}
\nonumber(A-B\ln x)x u'(x)+\frac{1}{2}u''(x)+\int_{\mathbb{R}
\setminus\{0\}}\frac{u(x+y)-u(x)-y
u'(x)I_{|y|<1}}{|y|^{1+\alpha}}\,dy&=&-1,
\quad  x\in D, \\
u(x)&=&0, \quad x\notin D.  \label{MET}
\end{eqnarray}
 The mean exit time $u(x)$ is the the time when the tumor (initially with density $x$) becomes not diagnoseable (when $x(t)$ exits from the left boundary point $x=1$),  or becomes   malignant
(when $x(t)$ exits from the right boundary point $x=b$) and thus may
need medical intervention.

The escape probability $p(x)$, through the right boundary point, is
the likelihood that the   tumor density $X(t)$, initially at value
$x$ in $D$,
  first   escapes from $D$ through $x= 1. 074\times10^8$.  The $p(x)$ satisfies the following equation
 \cite{Brannan99}
\begin{eqnarray} \label{EP}
Ap(x)&=&0,\quad x\in D=(1, 1. 074\times10^8), \nonumber\\
p(1)&=&0, \quad   p( 1. 074\times10^8)=1.
\end{eqnarray}

In the next section we simulate both the mean exit time $u(x)$ and
escape probability $p(x)$, when the stability parameter $\alpha$
varies in $(0, 2)$.

\section{Numerical experiments}

 We now discretize the equation \eqref{MET} for the mean exit time $u(x)$.
 The discretization of \eqref{EP} for the escape probability $p(x)$ is similar.

For the numerical scheme as the right boundary point of $D$ is too
huge, we need a transformation to change it to a smaller value.
Owing to the fact that $u(x)=0$ for $x\notin D$ we can take a
transformation $\tilde x=\ln x$, relevantly $\tilde y=\ln (x+y)-\ln
x$, nonetheless the integration for $y$ is on
$\mathbb{R}\setminus\{0\}$ ($u(x+y)=0$ for $x+y\leq0<1$). Omitting
the tilde we get
\begin{eqnarray}
\nonumber\lefteqn{(A-Bx)u'(x)+\frac{1}{2}e^{-2x}(u''(x)-u'(x))}\quad\quad\\
{}&+&\nonumber\frac{1}{e^{\alpha
x}}\int_{\mathbb{R}\setminus\{0\}}\frac{u(x+y)-u(x)-(e^y-1)u'(x)I_{\{|e^y-1|<e^{-x}\}}}{|e^y-1|^{1+\alpha}}e^y\,dy=-1,
\quad x\in D_1=(0, 18. 4921), \\
&&\quad \quad\quad\quad\quad\quad
\quad\quad\quad\quad\quad\quad\quad\quad\quad\quad\quad\quad\quad\quad\quad\quad\quad\quad
u(x)=0,\quad x\notin D_1.\label{(0,18.4921)}
\end{eqnarray}
where the third term in the integral is
$$\int_{\mathbb{R}\setminus\{0\}}\frac{(e^y-1)I_{\{|e^y-1|<e^{-x}\}}}{|e^y-1|^{1+\alpha}}e^y\,dy
=\int_{\mathbb{R}\setminus\{0\}}\frac{(e^y-1)I_{\{|e^y-1|<e^{-x}\}}}{|e^y-1|^{1+\alpha}}\,d(e^y-1)
=\int_{\mathbb{R}\setminus\{0\}}\frac{yI_{\{|y|<e^{-x}\}}}{|y|^{1+\alpha}}\,dy,$$
as
$\int_{\mathbb{R}\setminus\{0\}}\frac{yI_{\{|y|<\delta\}}}{|y|^{1+\alpha}}\,dy$
always vanishes for any $\delta>0$, so does the third term in the
integral for any fixed $x$, therefore the integral-differential
equation can be changed into
\begin{eqnarray}
\nonumber\lefteqn{\frac{1}{2}e^{-2x}u''(x)+[(A-Bx)-\frac{1}{2}e^{-2x}]u'(x)}\\
{}&+&\nonumber\frac{1}{e^{\alpha
x}}\int_{\mathbb{R}\setminus\{0\}}\frac{u(x+y)-u(x)-(e^y-1)u'(x)I_{\{|e^y-1|<\delta\}}}{|e^y-1|^{1+\alpha}}e^y\,dy=-1,
\quad x\in D_1,\\
&&\quad\quad\quad\quad\quad\quad\quad\quad\quad
\quad\quad\quad\quad\quad\quad\quad\quad\quad\quad\quad\quad\quad\quad
u(x)=0, \quad x\notin D_1.
\end{eqnarray}
We then take another transformation $\tilde x=\frac{x}{b}-1$
meanwhile $\tilde y=\frac{y}{b}$ where $b=18. 4921/2$ to change
$D_1=(0, 18. 4921)$ into $(-1, 1)$, we obtain
\begin{eqnarray}
\nonumber\lefteqn{\frac{1}{2b^2}e^{-2b(x+1)}u''(x)+\big[A-B
b(x+1)-\frac{1}{2}e^{-2b(x+1)}\big]\frac{1}{b}u'(x)}\\
{}&+&\nonumber\frac{1}{e^{\alpha b
(x+1)}}\int_{\mathbb{R}\setminus\{0\}}\frac{b u(x+y)-b
u(x)-(e^{by}-1)u'(x)I_{\{|e^{by}-1|<\delta\}}}{|e^{by}-1|^{1+\alpha}}e^{by}\,dy=-1,
\quad x\in D_2=(-1, 1),\\
&&{}\quad\quad\,\,\,\quad\quad\quad\quad\quad\quad\quad\quad\quad\quad\quad\quad\quad\quad\quad\quad\quad\quad\quad\quad\quad\quad\quad\quad
u(x)=0, \quad x\notin D_2.
\end{eqnarray}
The numerical scheme here follows \cite{Li11}. A different scheme
was also proposed in \cite{Chen11}. We write
$\int_{\mathbb{R}}=\int_{-\infty}^{-1-x}+\int_{-1-x}^{1-x}+\int_{1-x}^\infty$
and take $\delta=\min{\big\{e^{b(1-x)}-1, 1-e^{-b(1+x)}\big\}}$ to
get
\begin{eqnarray}
\nonumber\lefteqn{\frac{1}{2b^2}e^{-2b(x+1)}u''(x)+\Big[A-B
b(x+1)-\frac{1}{2}e^{-2b(x+1)}\Big]\frac{1}{b}u'(x)}\\{}&+&\frac{u(x)}{\alpha
e^{\alpha b
(x+1)}}\Big[1-\frac{1}{(e^{b(1-x)}-1)^\alpha}-\frac{1}{(1-e^{-b(1+x)})^\alpha}\Big]\nonumber\\
{}&+&\frac{1}{ e^{\alpha b (x+1)}}\int_{-1-x}^{1-x}\frac{b u(x+y)-b
u(x)-(e^{by}-1)u'(x)I_{\{|e^{by}-1|<\delta\}}}{|e^{by}-1|^{1+\alpha}}e^{by}\,dy=-1,\quad
x\in
D_3=(-1, 1)\nonumber\\
\end{eqnarray}
furthermore, for $x\geq0$
\begin{eqnarray}
\nonumber\lefteqn{\frac{1}{2b^2}e^{-2b(x+1)}u''(x)+\Big[A-B
b(x+1)-\frac{1}{2}e^{-2b(x+1)}\Big]\frac{1}{b}u'(x)}\\{}&+&\frac{u(x)}{\alpha
e^{\alpha
b(x+1)}}\Big[1-\frac{1}{(e^{b(1-x)}-1)^\alpha}-\frac{1}{(1-e^{-b(1+x)})^\alpha}\Big]+
\frac{1}{e^{\alpha b (x+1)}}\int_{-1-x}^{-1+x}\frac{u(x+y)-
u(x)}{|e^{by}-1|^{1+\alpha}}be^{by}\,dy\nonumber\\{}&+&\frac{1}{e^{\alpha
b (x+1)}}\int_{-1+x}^{1-x}\frac{b u(x+y)-b
u(x)-(e^{by}-1)u'(x)}{|e^{by}-1|^{1+\alpha}}e^{by}\,dy=-1,\quad x\in
D_3=(-1, 1)\label{u(x)plus}
\end{eqnarray}
and for $x\leq0$
\begin{eqnarray}
\nonumber\lefteqn{\frac{1}{2b^2}e^{-2b(x+1)}u''(x)+\Big[A-B
b(x+1)-\frac{1}{2}e^{-2b(x+1)}\Big]\frac{1}{b}u'(x)}\\{}&+&\frac{u(x)}{\alpha
e^{\alpha
b(x+1)}}\Big[1-\frac{1}{(e^{b(1-x)}-1)^\alpha}-\frac{1}{(1-e^{-b(1+x)})^\alpha}\Big]+
\frac{1}{e^{\alpha b (x+1)}}\int_{1+x}^{1-x}\frac{u(x+y)-
u(x)}{|e^{by}-1|^{1+\alpha}}be^{by}\,dy\nonumber\\{}&+&\frac{1}{e^{\alpha
b (x+1)}}\int_{-1-x}^{1+x}\frac{b u(x+y)-b
u(x)-(e^{by}-1)u'(x)}{|e^{by}-1|^{1+\alpha}}e^{by}\,dy=-1,\quad x\in
D_3=(-1, 1)\label{u(x)minus}
\end{eqnarray}
For simplicity denote $f(x)=A-B b(x+1)-\frac{1}{2}e^{-2b(x+1)}$ and
$g(x)=e^{-2b(x+1)}$. For numerical schemes generally we divide the
interval $[-2,2]$ into $4J$ subintervals and define $x_j=jh$ for
$-2J\leq j \leq 2J$ integer, where $h=1/J$. We denote the numerical
solution of $u$ at $x_j$ by $U_j$. We can discretize the two
integral-differential equations above using central difference for
derivatives and ``punched-hole'' trapezoidal rule
\begin{eqnarray}
\lefteqn{\frac{1}{2b^2}g(x_j)\frac{U_{j+1}-2U_j+U_{j-1}}{h^2}+\frac{1}{b}f(x_j)\frac{U_{j+1}-U_{j-1}}{2h}}\nonumber\\{}&+&\frac{1}{\alpha
e^{\alpha
b(x_j+1)}}\Big[1-\frac{1}{(e^{b(1-x_j)}-1)^\alpha}-\frac{1}{(1-e^{-b(1+x_j)})^\alpha}\Big]U_j+\frac{h}{e^{\alpha
b(x_j+1)}}\sum_{k=-J-j}^{-J+j}{''}\,\,
\frac{be^{by_k}(U_{j+k}-U_j)}{|e^{by_k}-1|^{1+\alpha}}\nonumber\\{}
&+&\frac{h}{e^{\alpha
b(x_j+1)}}\sum_{k=-J+j,\,\,k\neq0}^{J-j}{''}\,\,
\frac{be^{by_k}(U_{j+k}-U_j)-e^{by_k}(e^{by_k}-1)
(U_{j+1}-U_{j-1})/2h}{|e^{by_k}-1|^{1+\alpha}}=-1,
\end{eqnarray}
where $j=0,1,2,\cdots,J-1$. Here the summation symbol $\sum{''}$
means the two end terms are multiplied by $1/2$ corresponding to
trapezoidformula.
\begin{eqnarray}
\lefteqn{\frac{1}{2b^2}g(x_j)\frac{U_{j+1}-2U_j+U_{j-1}}{h^2}+\frac{1}{b}f(x_j)\frac{U_{j+1}-U_{j-1}}{2h}}\nonumber\\{}
&+&\frac{1}{\alpha e^{\alpha
b(x_j+1)}}\Big[1-\frac{1}{(e^{b(1-x_j)}-1)^\alpha}-\frac{1}{(1-e^{-b(1+x_j)})^\alpha}\Big]U_j
+\frac{h}{e^{\alpha b(x_j+1)}}\sum_{k=J+j}^{J-j}{''}\,\,
\frac{be^{by_k}(U_{j+k}-U_j)}{|e^{by_k}-1|^{1+\alpha}}\nonumber\\{}
&+&\frac{h}{e^{\alpha b(x_j+1)}}\sum_{k=J+j,\,\,k\neq0}^{-J-j}{''}\,
\,\frac{be^{by_k}(U_{j+k}-U_j)-e^{by_k}(e^{by_k}-1)
(U_{j+1}-U_{j-1})/2h}{|e^{by_k}-1|^{1+\alpha}}=-1,\quad\quad\quad\quad\quad
\end{eqnarray}
where $j=-(J-1),-(J-2),\cdots,-2,-1$. The boundary conditions
require that the values of $U_j$ vanish if the index $|j|\geq J$.
But here for the specific problem for first-order derivative we use
forward difference or backward difference.

\medskip

Similarly, we can discretize the equation \eqref{EP} for the escape
probability $p(x)$.

\medskip

For the viewing convenience, the following figures for the mean exit
time $u(x)$  and escape probability $p(x)$ are both plotted for
$x\in(0, 18.4921)$,   instead of the original huge range $(1,
1.074\times10^8)$.
\medskip

\begin{figure}[h]
\psfig{file=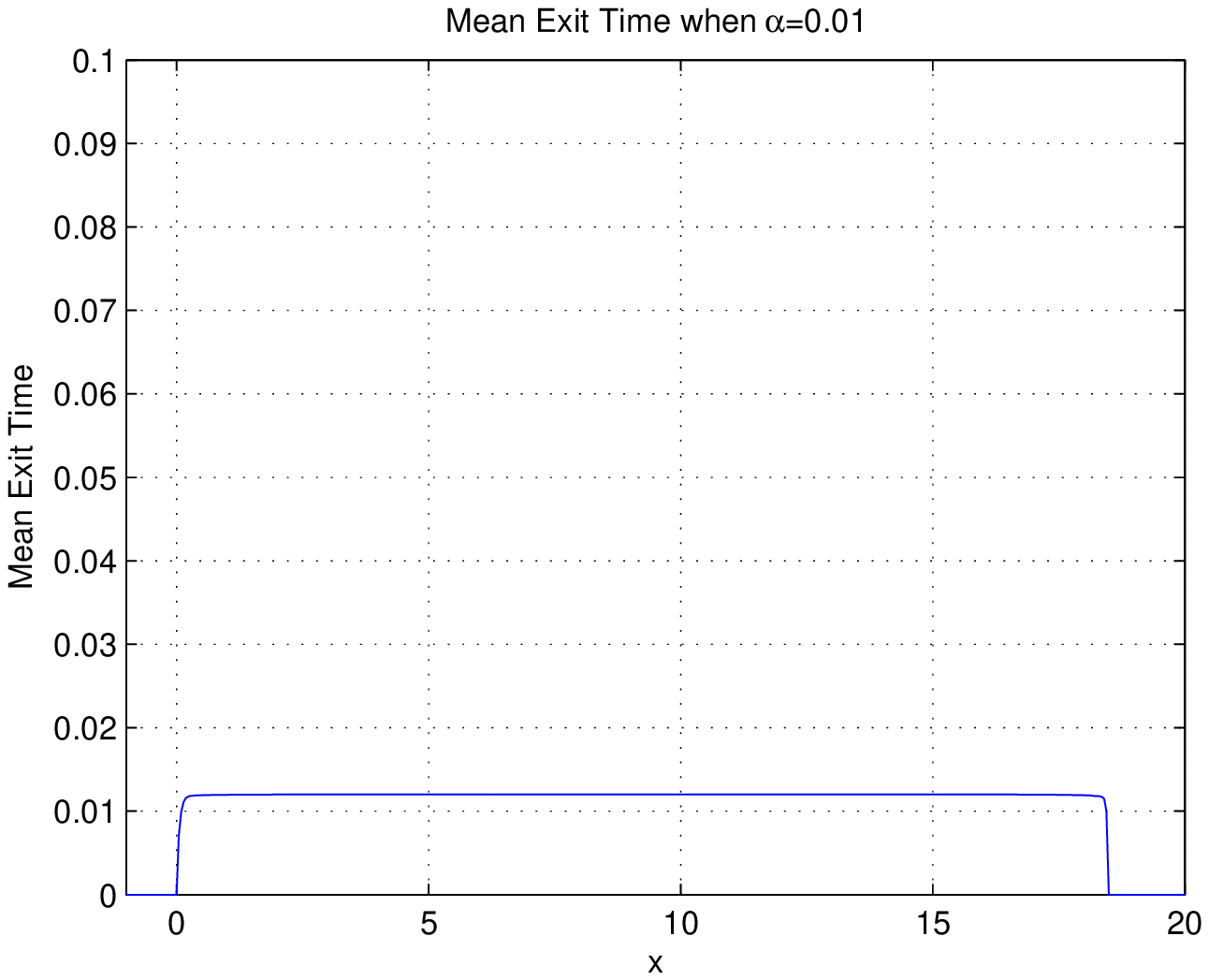} \psfig{file=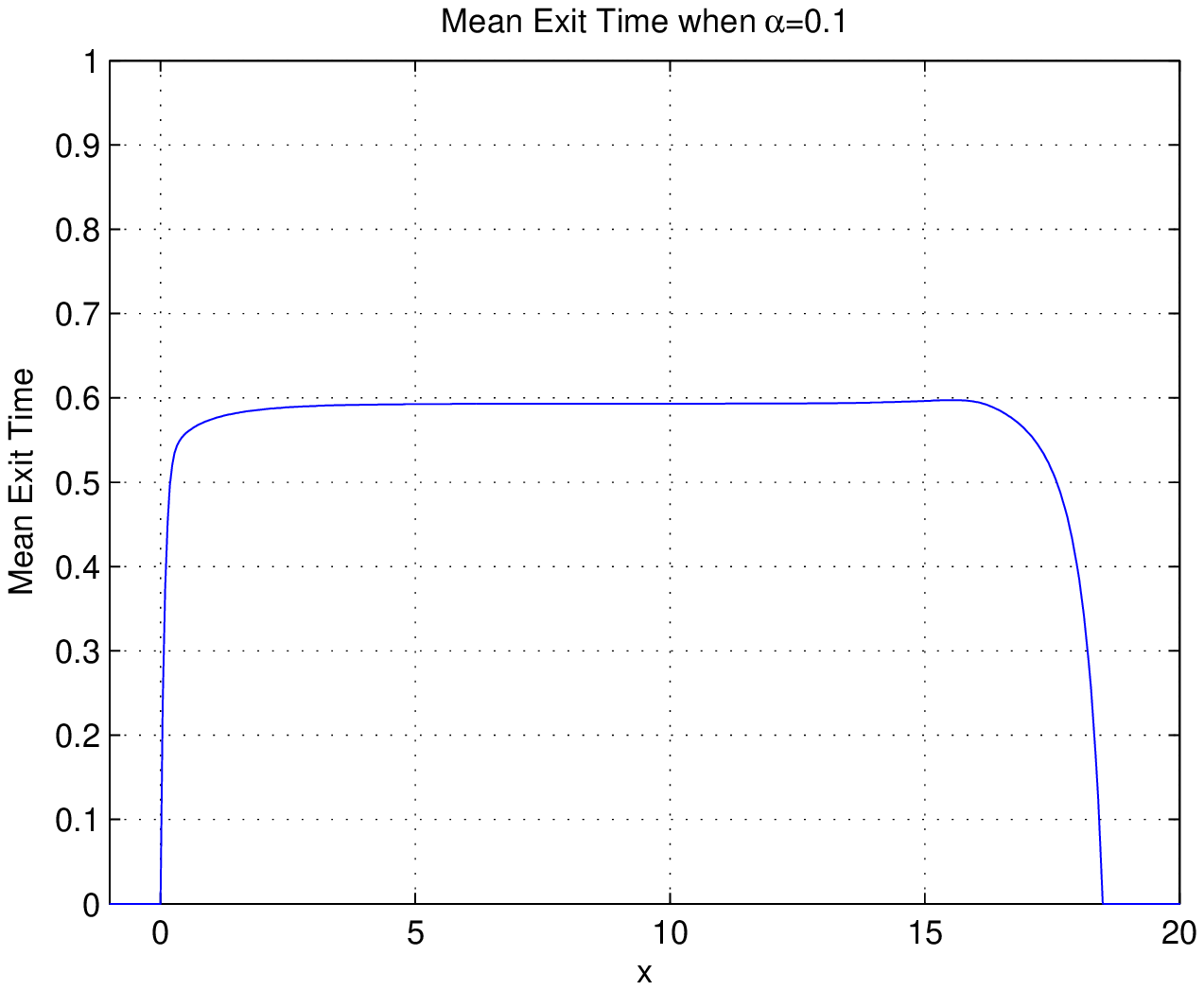}
\caption{Plot of $u(x)$ given by (\ref{(0,18.4921)}) with
$\alpha=0.01$ and $\alpha=0.1$. } \label{MET1}
\end{figure}

\begin{figure}[h]
\psfig{file=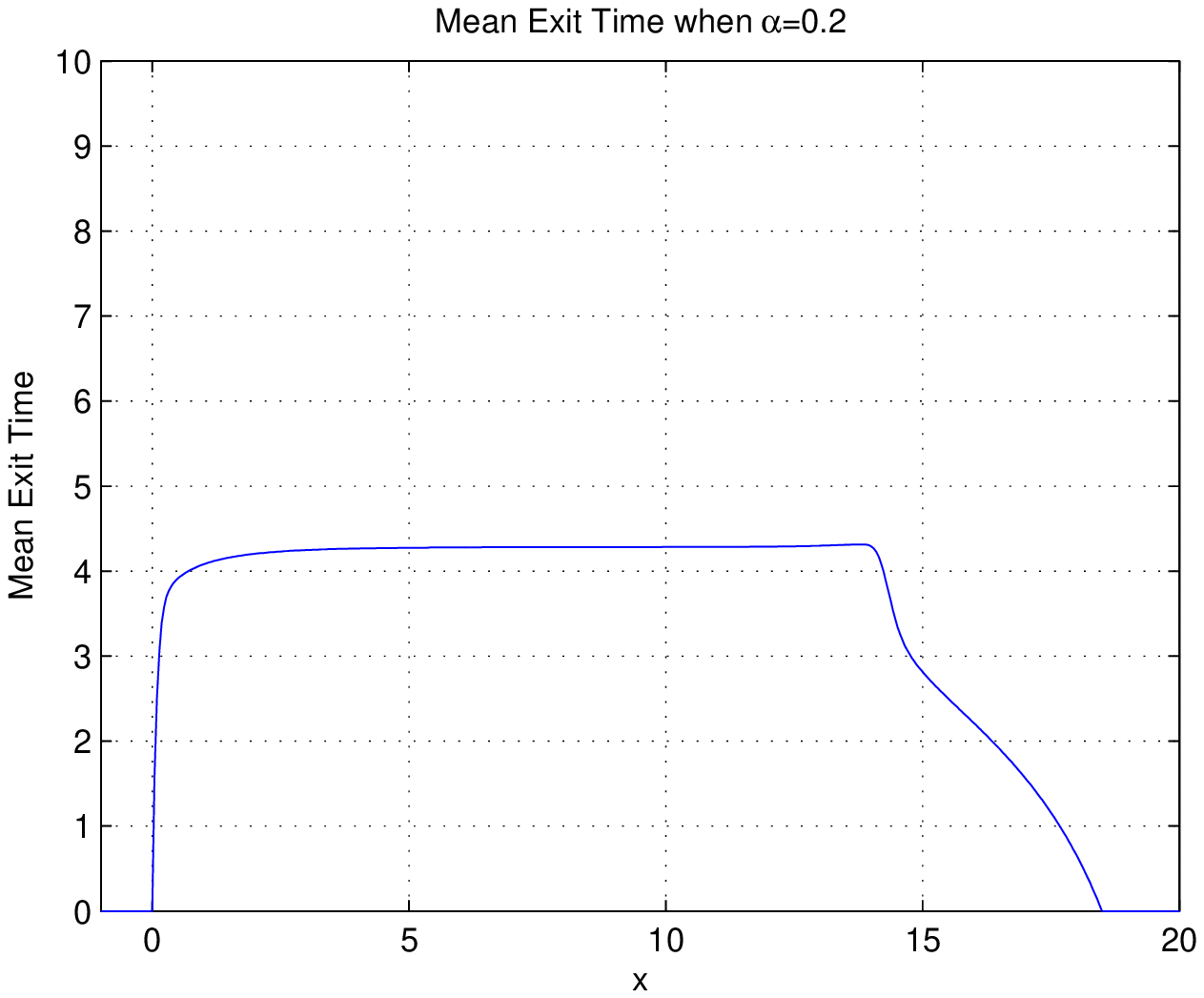} \psfig{file=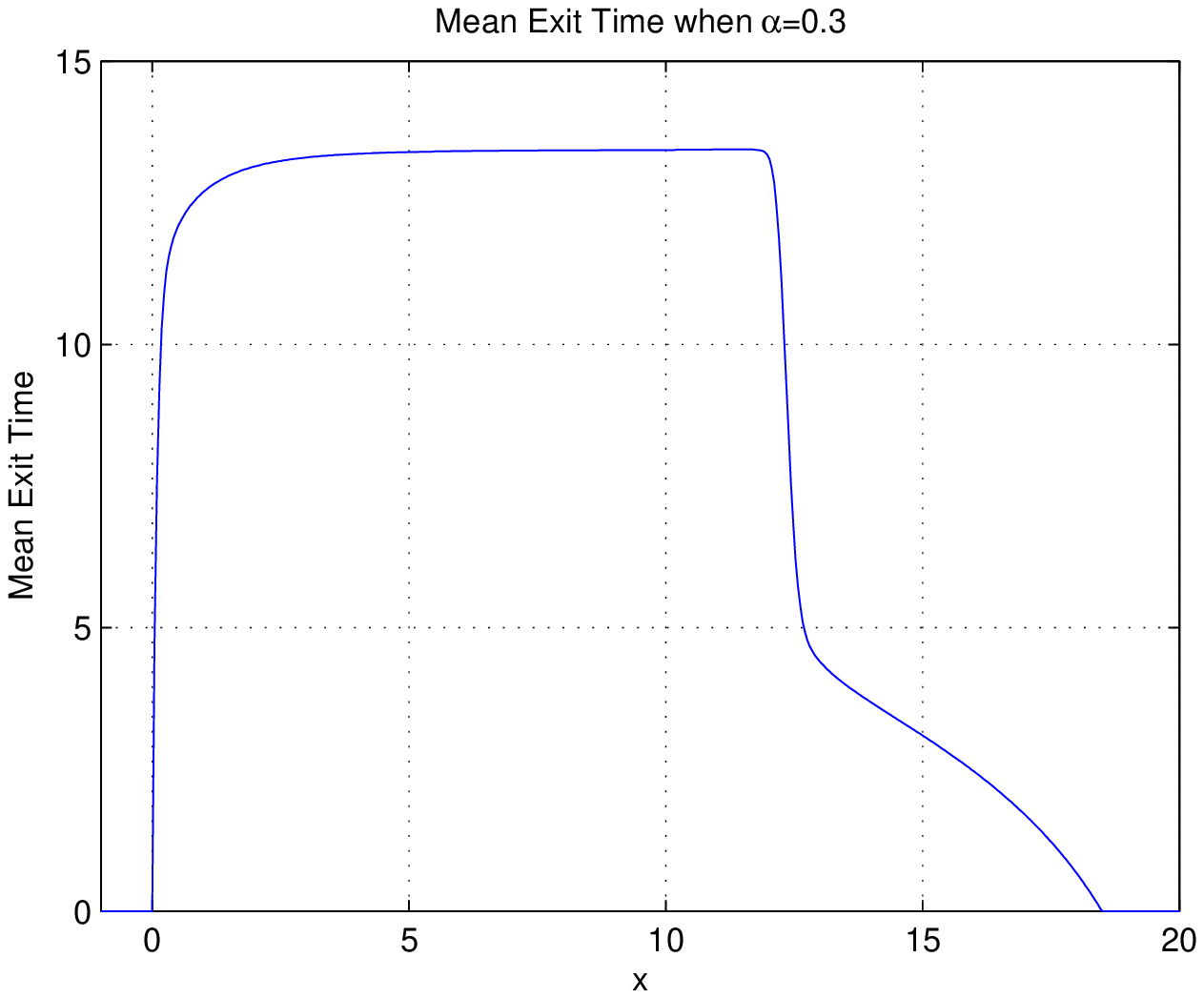}
\caption{Plot of $u(x)$ given by (\ref{(0,18.4921)}) with
$\alpha=0.2$ and $\alpha=0.3$. }\label{MET2}
\end{figure}

\begin{figure}[h]
\psfig{file=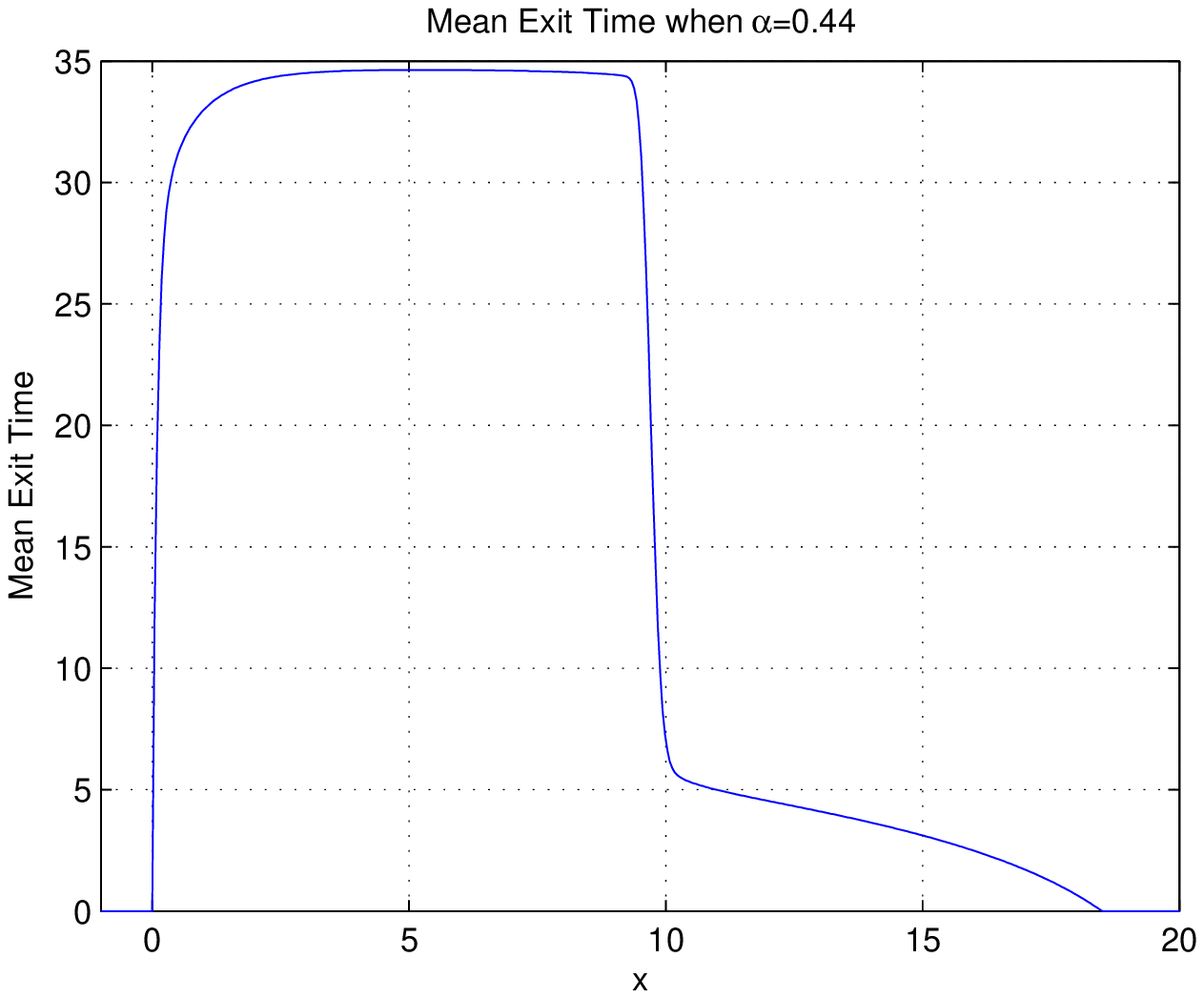} \psfig{file=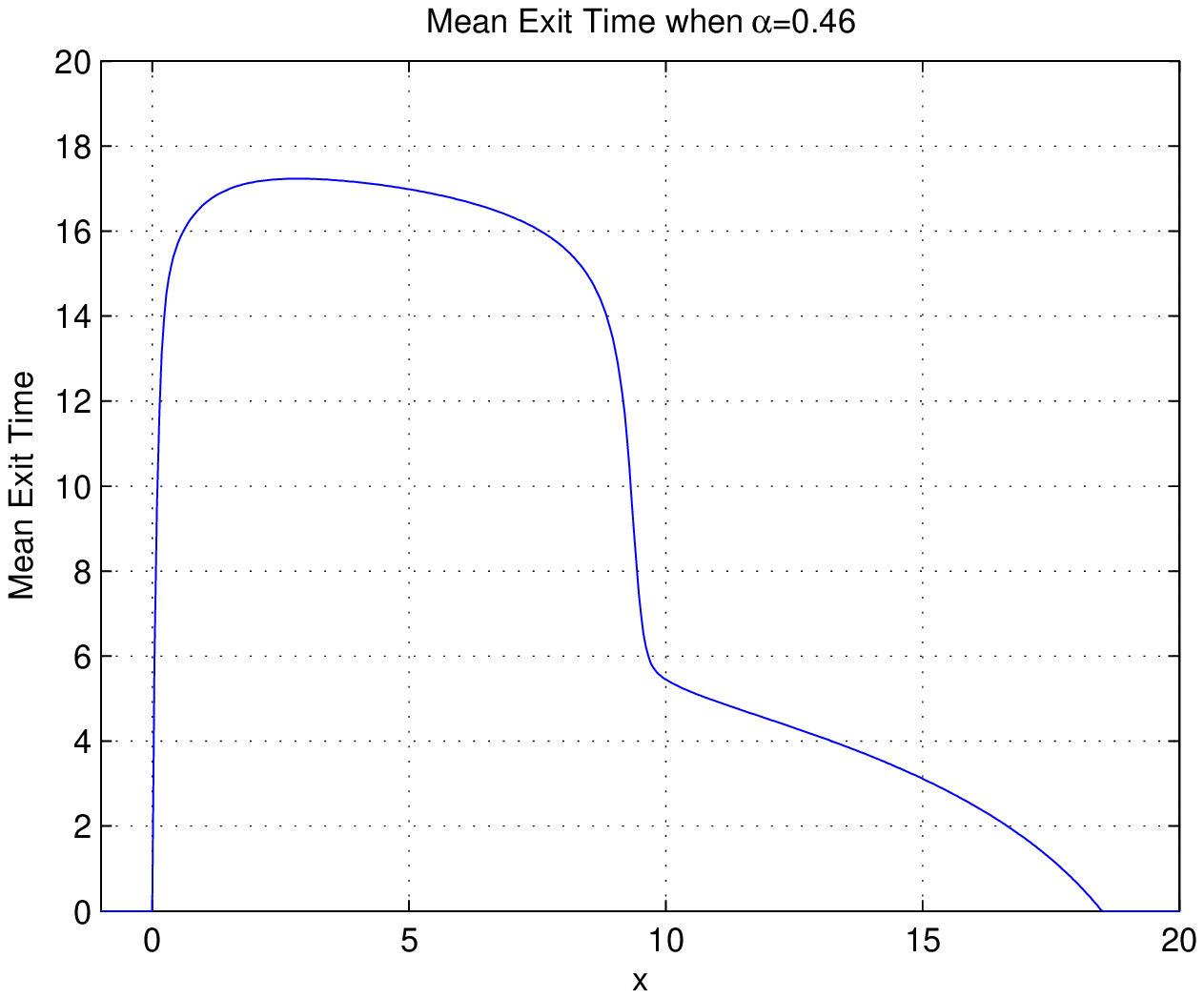}
\caption{Plot of $u(x)$ given by (\ref{(0,18.4921)}) with
$\alpha=0.44$ and $\alpha=0.46$.}\label{MET3}
\end{figure}

\begin{figure}[h]
\psfig{file=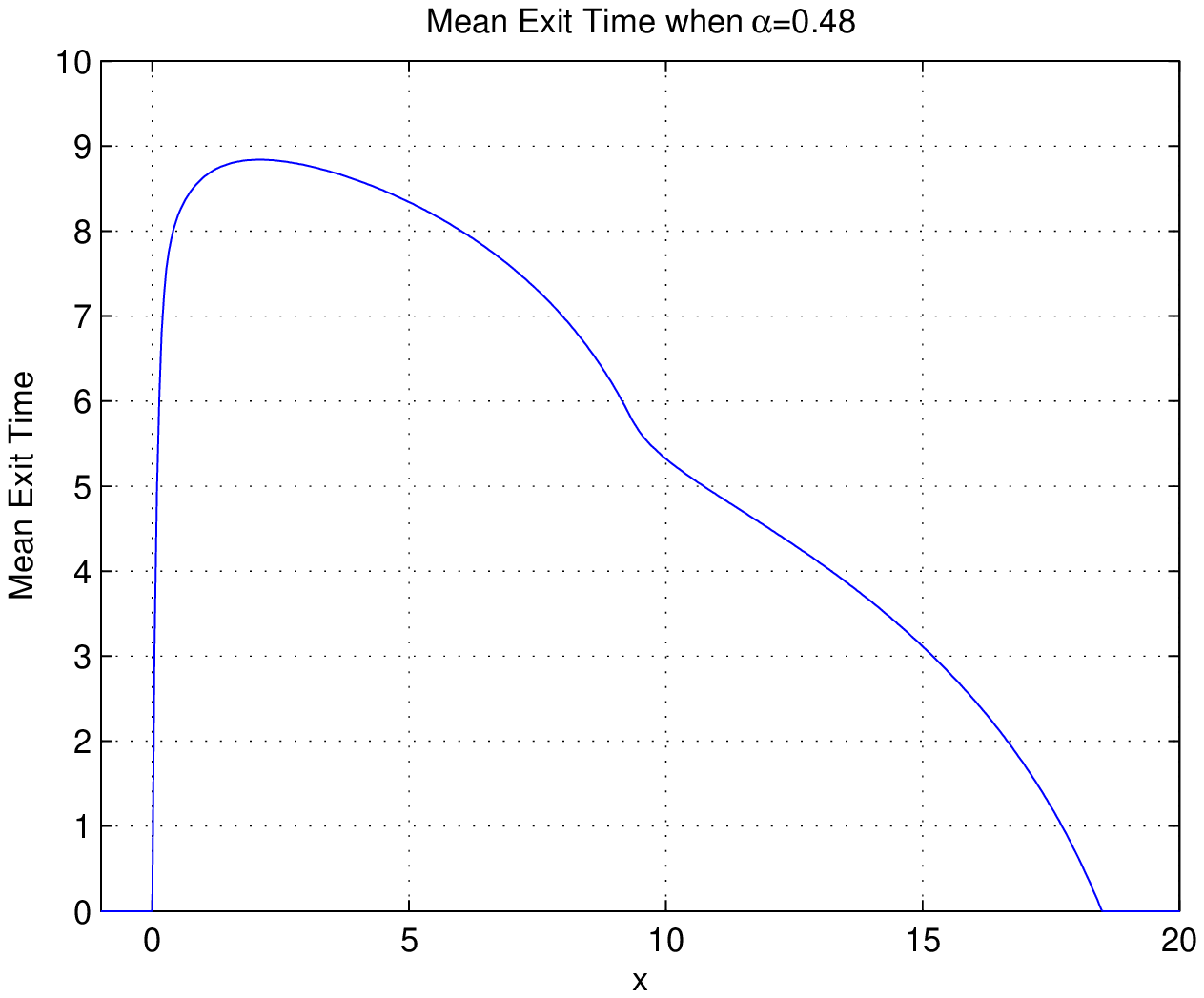} \psfig{file=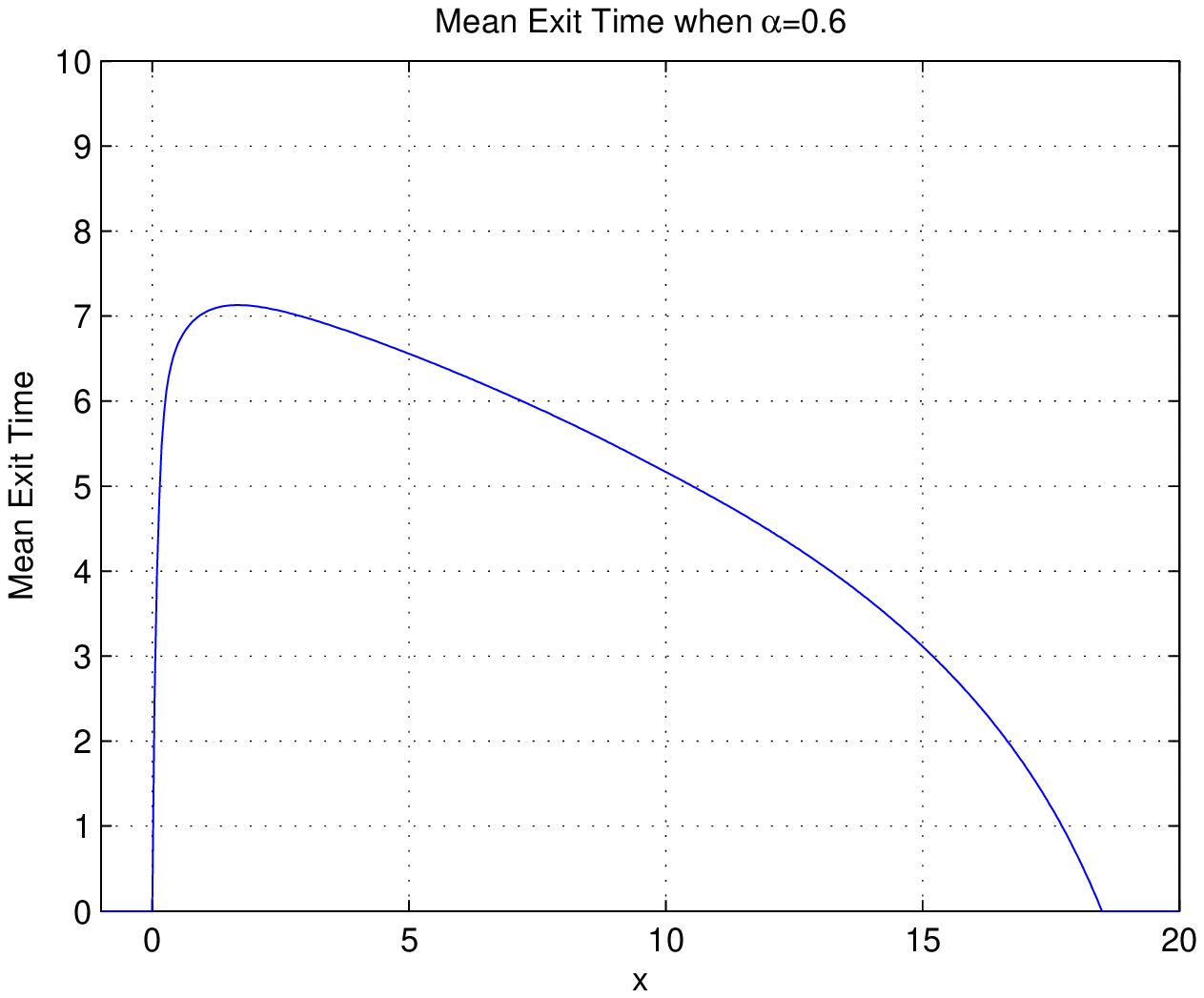}
\caption{Plot of $u(x)$ given by (\ref{(0,18.4921)}) with
$\alpha=0.48$ and $\alpha=0.6$. }\label{MET4}
\end{figure}

\begin{figure}[h]
\psfig{file=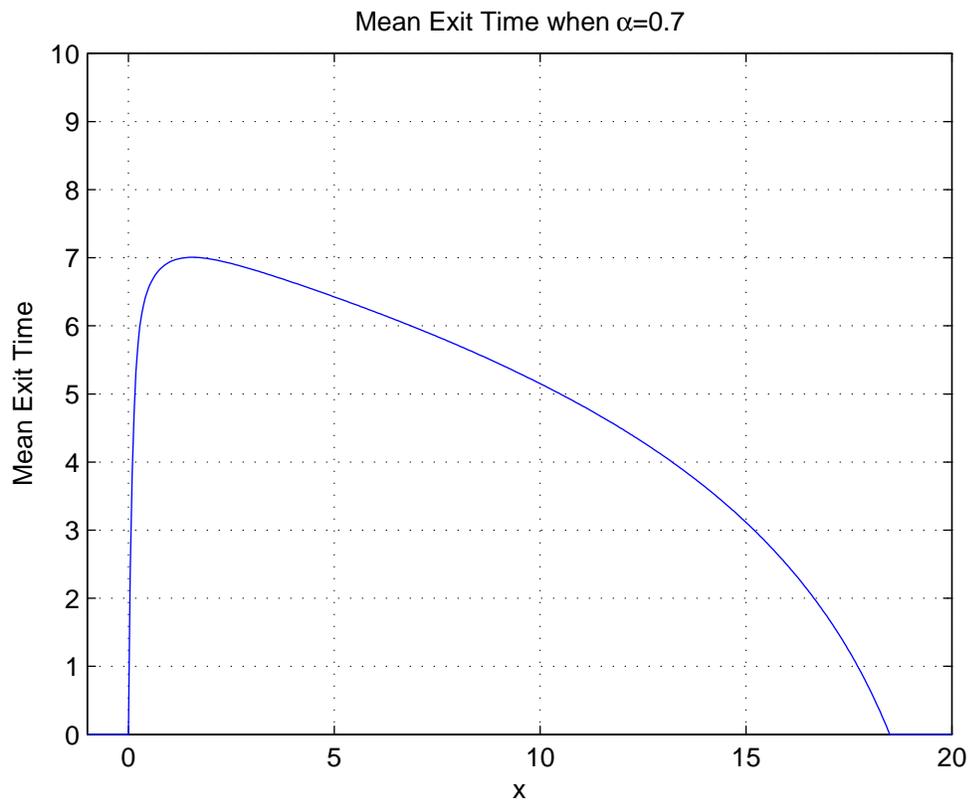} \psfig{file=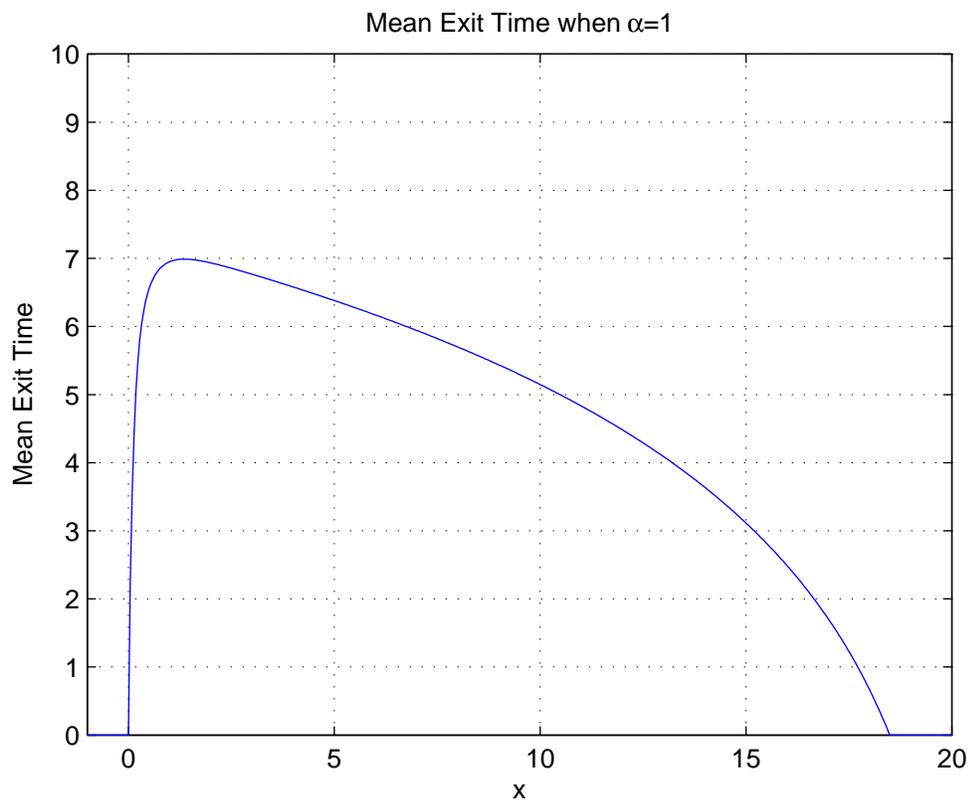}
\caption{Plot of $u(x, t)$ given by (\ref{(0,18.4921)}) with
$\alpha=0.7$ and $\alpha=1$. }\label{MET5}
\end{figure}

\begin{figure}[h]
\psfig{file=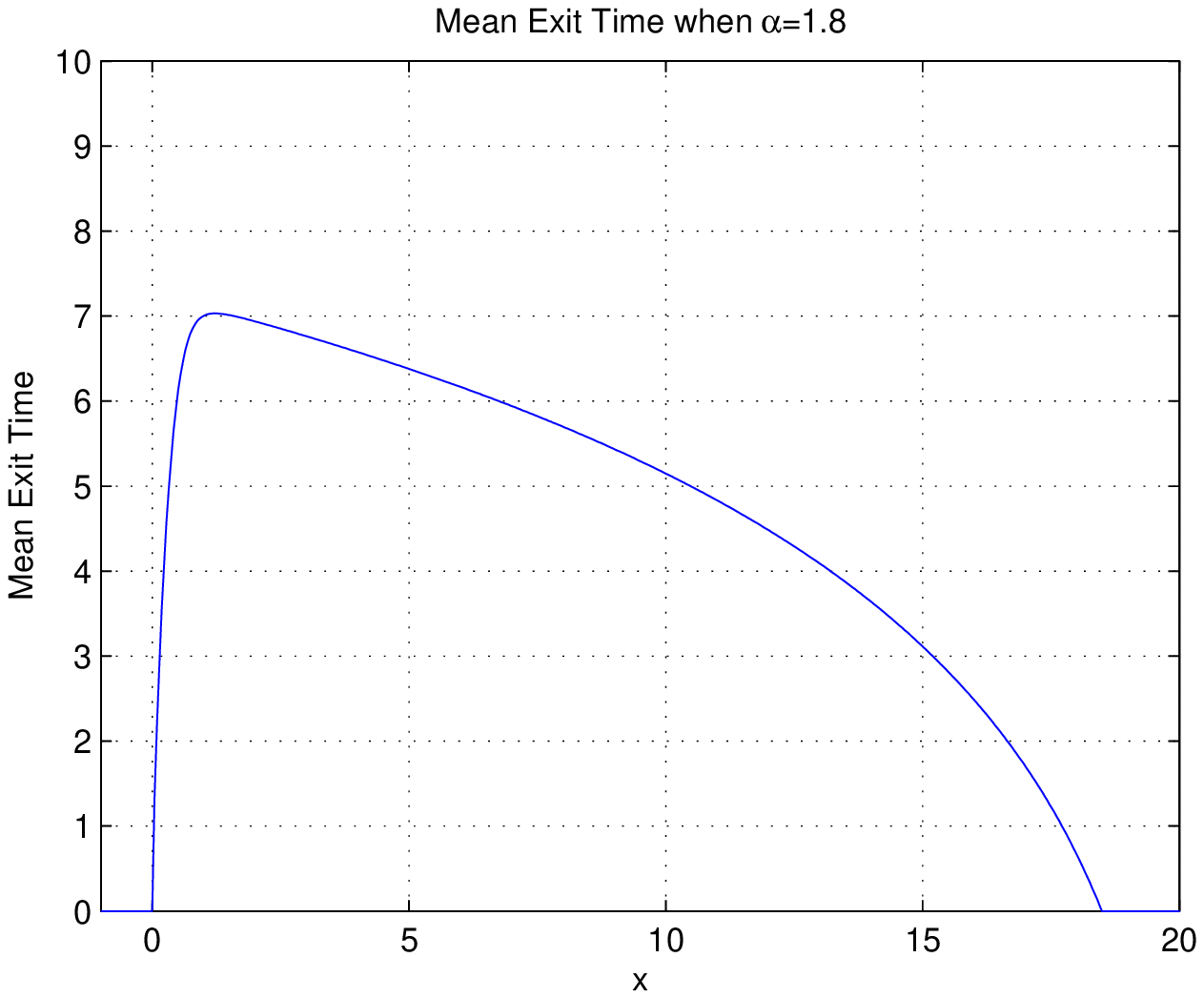} \psfig{file=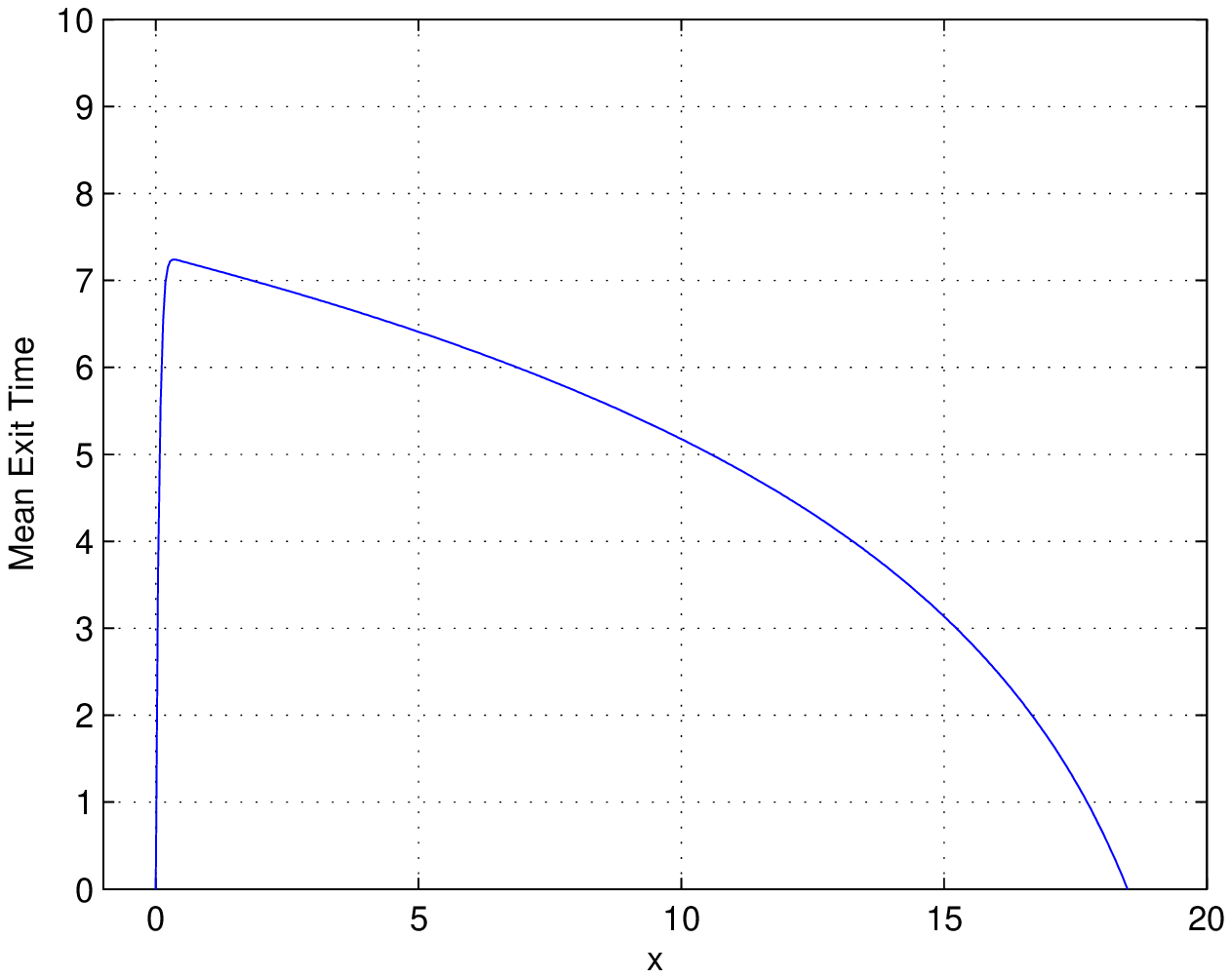} \caption{Plot of
$u(x, t)$ given by (\ref{(0,18.4921)}) with $\alpha=1.8$ and  the
Brownian motion case. }   \label{MET6}
\end{figure}

Figures \ref{MET1}--\ref{MET6} show the mean exit time $u(x)$ for
various $\alpha$ values in $(0, 2)$.

We observe that when $\alpha$ is small, the mean exit time $u(x)$ is
small and   is roughly  constant, i.e.,  tumors at all densities are
equally unlikely to become non-diagnosable or
 malignant. As $\alpha$
becomes bigger the whole height will grow with the right part grows
up slowly while the left part increases rapidly meanwhile the width
of the left part narrows down. In fact, at about $\alpha=0.2$ the
height of the two part begins to have a separation. When $\alpha
\approx 0.44$ the separation becomes quite obvious and the left part
(for density between $0$ and $10$) reaches its highest value, which
indicates a bifurcation. This means that for a tumor with density
between $0$ and $10$, it takes longest time to become either
non-diagnosable or  malignant.

For $\alpha$ between $0.44 $ and $0.48$, $u(x)$ keeps similar shape
as for $\alpha=0.44$ but the peak diminishes. For $\alpha$ between
$0.48 $ and $0.6$, $u(x)$ has the similar shape although the peak
value decreases. Moreover, for $\alpha$ between $0.6 $ and $2$,
$u(x)$ keeps almost the same shape; namely,   for the Gaussian
Brownian motion case, the mean exit time is similar to the large
values of $\alpha \in (0.6, 2)$. The dramatic impact of non-Gaussian
L\'evy noise occurs for $\alpha \in (0, 0.46)$.

\begin{figure}[h]
\begin{center}
 \psfig{file=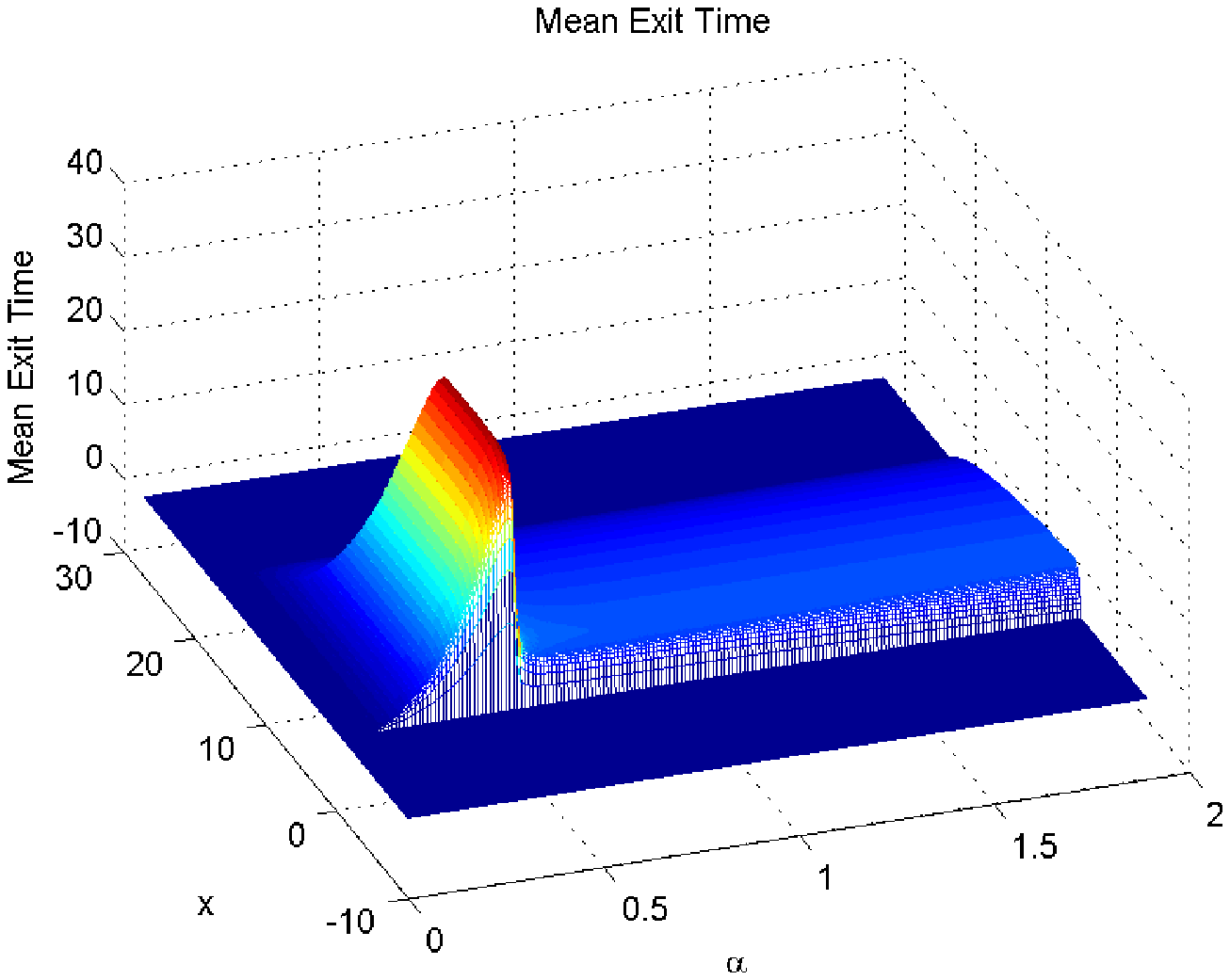}
\end{center}
\caption{Mean exit time $u$ v.s. $\alpha$ and $x$. }
\label{thr-dim}
\end{figure}

Figure \ref{thr-dim} shows a three dimensional plot of the above
discussed situations.

%%%%%%%%%%%%%%%%%%%%%%%%%%%%%%%%%%%%%%%%%%%%%%%%%%%%%%%%%%%%%%%%
%%%%%%%%%%%%%%%%%%%%%%%%%%%%%%%%%%

\begin{figure}[h]
\psfig{file=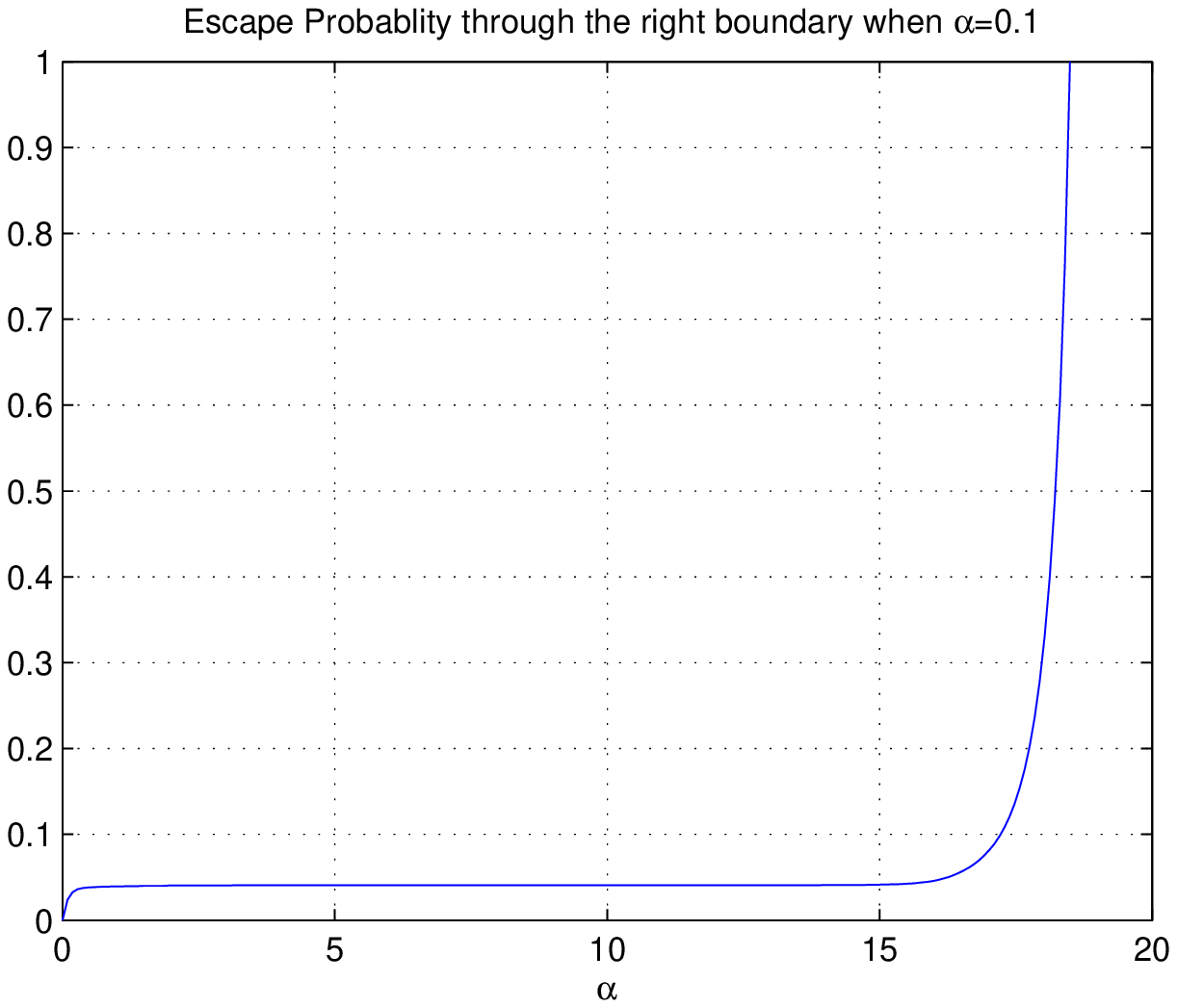}
\psfig{file=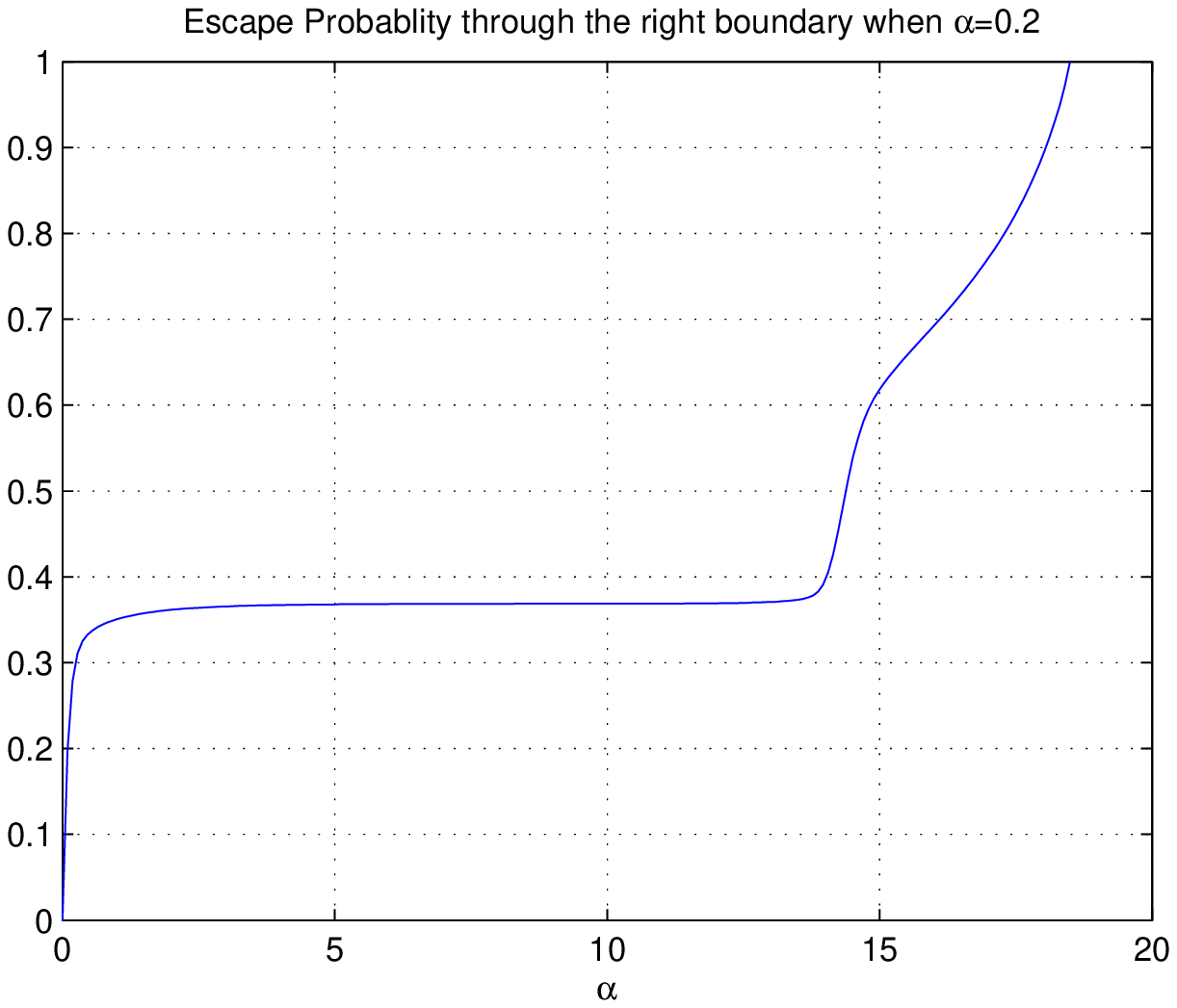}
\caption{Escape
probability $p(x)$ through the right boundary of $(0, 18.4921)$:
$\alpha=0.1$ and   $\alpha=0.2$. }  \label{escape1}
\end{figure}

\begin{figure}[h]
\psfig{file=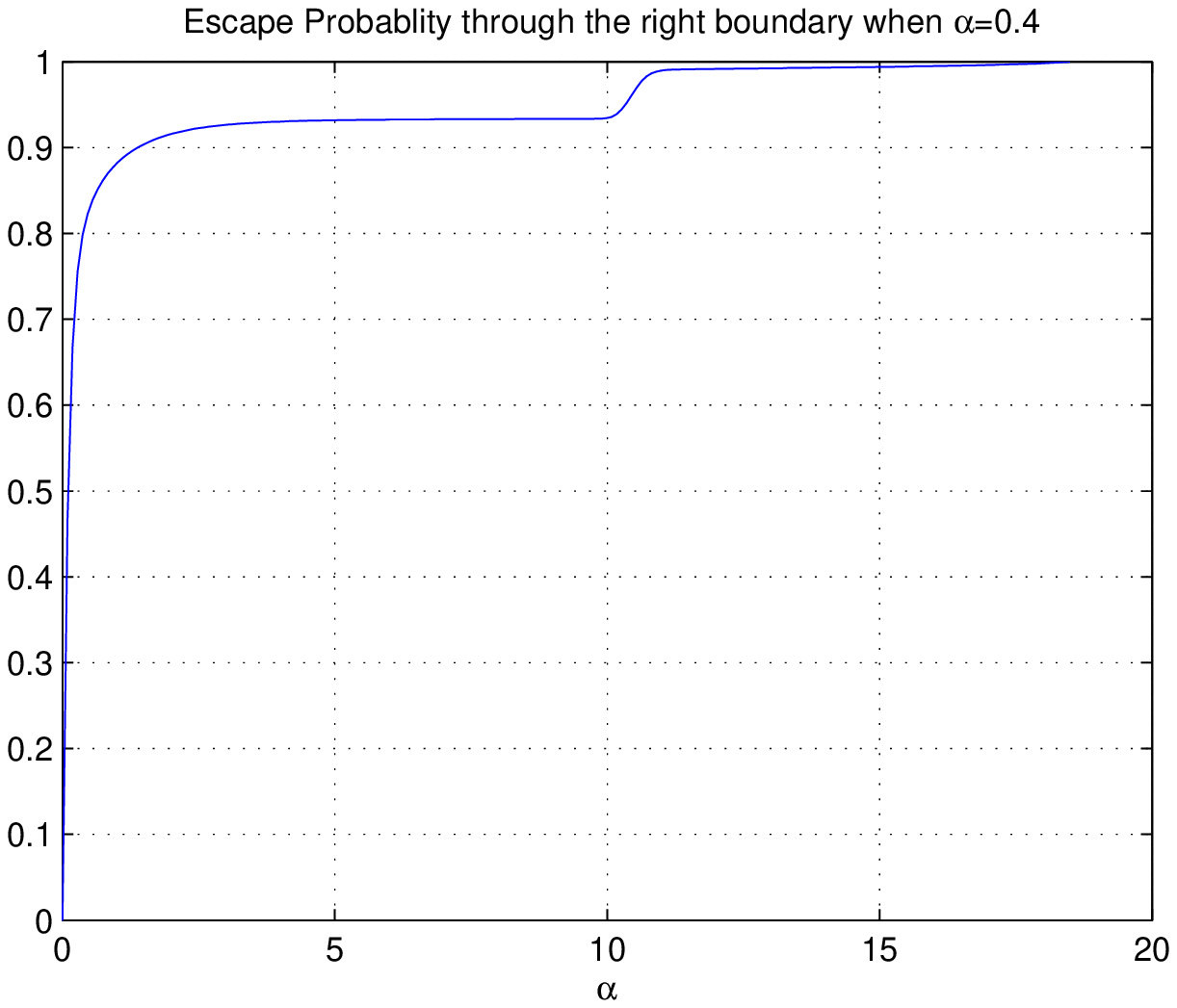}
\psfig{file=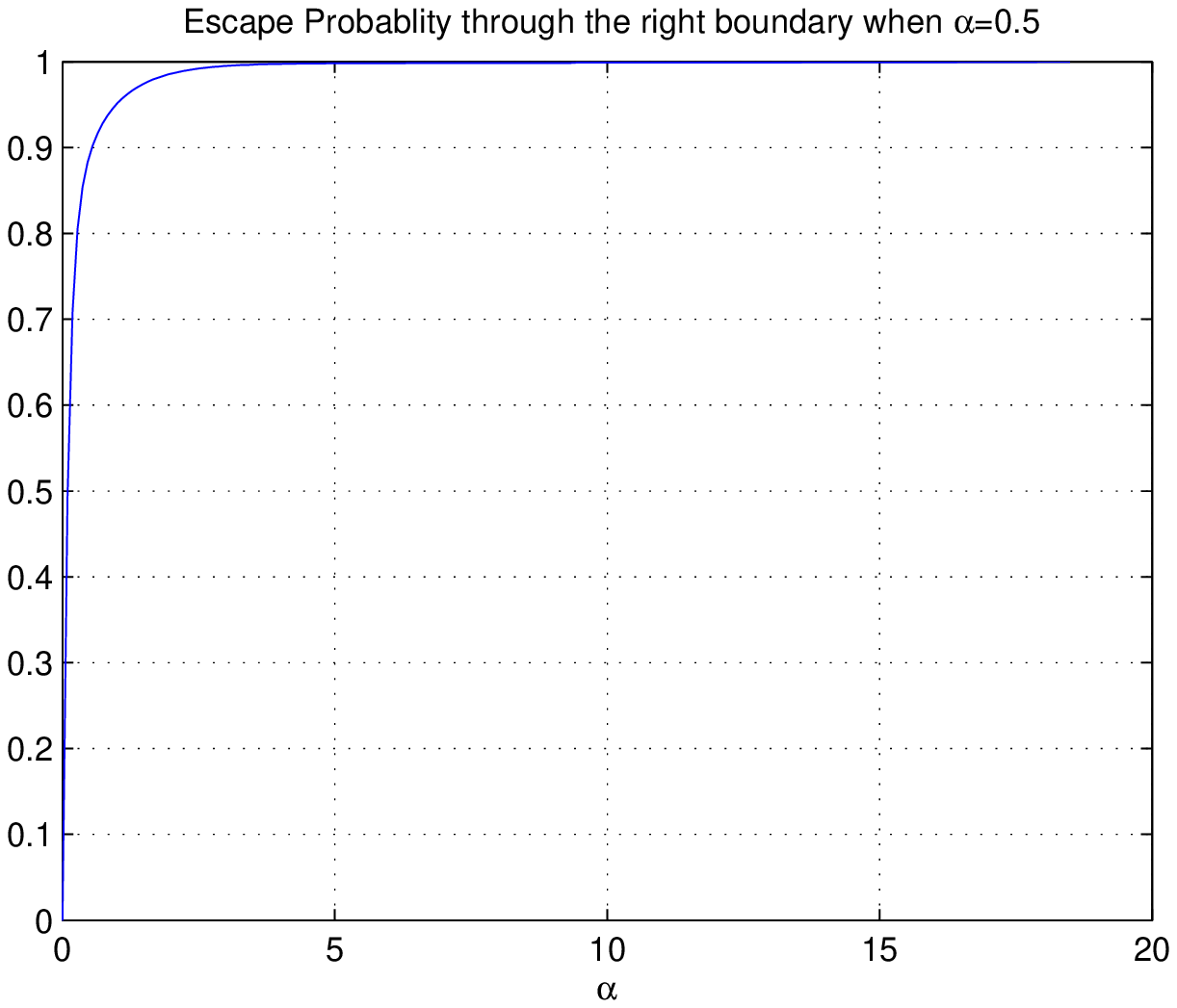}
\caption{Escape
probability $p(x)$ through the right boundary of $(0, 18.4921)$:
$\alpha=0.4$ and   $\alpha=0.5$. } \label{escape2}
\end{figure}

In the  above simulations, the  mean exit time $u(x)$ tells us the
time for a tumor to become either non-diagnosable  (exit from the
left boundary point) or become malignant (exit from the  right
boundary point). To distinguish these two situations, especially
examine the situation for a tumor to become malignant, let us now
consider
  the escape probability, $p(x)$, through the right end point. It indicates the likelihood the the tumor
   is progressing from benign to  malignant.   We observe that a bifurcation occurs for $\alpha \approx 0.4$ when a tumor at any benign density will highly likely becomes malignant, while for $\alpha < 0.4$, only tumors with high density  (near $x=20$) are likely to become malignant; see Figures \ref{escape1} and \ref{escape2}.

%%%%%%%%%%%%%%%%%%%%%%%%%%%%%
%%%%%%%%%%%%%%%%%%%%%%%%

\medskip

\noindent {\bf Acknowledgements.} We thank Xiaofan Li and Ting Gao
for help with numerical schemes, and Xu Sun for helpful discussions
on the model.

%\end{multicols}
\end{document}